\documentclass[11pt,a4paper]{amsart}

\usepackage{amsmath}
\usepackage{amssymb}
\usepackage{amsthm}
\usepackage{enumerate}
\usepackage{mathrsfs} 
\usepackage[T1]{fontenc}
\usepackage[utf8]{inputenc}
\usepackage{xcolor}
\usepackage[unicode]{hyperref}
\hypersetup{
colorlinks={true},
 linkcolor={red},
 citecolor={green!60!black},
  urlcolor={blue},
}

\theoremstyle{plain}
 \newtheorem{thm}{Theorem}
 \newtheorem*{thm*}{Theorem}
 
 \newtheorem*{cor*}{Corollary}
 \newtheorem{lem}[thm]{Lemma}
 \newtheorem*{lem*}{Lemma}
 \newtheorem{prop}[thm]{Proposition}
 \newtheorem*{prop*}{Proposition}
 
 \newtheorem*{obs*}{Observation}
 \newtheorem*{tw*}{Twierdzenie}
\theoremstyle{definition}
 \newtheorem{rem}[thm]{Remark}
 \newtheorem*{rem*}{Remark}
 \newtheorem{exmp}[thm]{Example}
 \newtheorem*{exmp*}{Example}
 \newtheorem*{ack*}{Acknowledgment}
 
 \newtheorem*{defn*}{Definition}

\newcommand{\N}{\mathbb{N}}

\newcommand{\R}{\mathbb{R}}

\newcommand{\itemref}[1]{\eqref{#1}} 

\renewcommand{\C}{\mathcal{C}}

\newcommand{\dx}[1][x]{\,\mathrm{d}#1}

\newcommand{\info}[1]{%
 \makeatletter
 \AtBeginDocument{\def\@serieslogo{\vbox{\noindent\Small #1}}\copyrightinfo{}{}}
 \makeatother}
\newcommand{\infowithcopyright}[3]{%
 \makeatletter
 \AtBeginDocument{\def\@serieslogo{\vbox{\noindent\Small #1}}\copyrightinfo{#2}{#3}}
 \makeatother}

\newcommand{\st}{such that}
\newcommand{\cf}{cf.}
\newcommand{\eg}{e.g.}
\newcommand{\ie}{i.e.}

\usepackage{fourier}

\renewcommand{\le}{\leqslant}

\renewcommand{\ge}{\geqslant}

\newcommand{\E}{\mathbb{E}\,}
\newcommand{\lcx}{\leqslant_{\normalfont{\textrm{cx}}}}

\newcommand{\xim}[1][(m)]{\xi_{#1}}
\newcommand{\gm}{\gamma_{(m)}}

\numberwithin{equation}{section}

\info{Manuscript submitted to J. Math. Anal. Appl.}

\theoremstyle{plain}
\newtheorem*{pr}{Problem}
\newtheorem*{Ohlin}{Ohlin's Lemma}
\newtheorem*{LevSte}{Levin--Ste\v{c}kin Theorem}

\begin{document}

\date{\today}
	\title{A solution to the problem of Ra\c{s}a connected with Bernstein polynomials}
	\author[J. Mrowiec]{Jacek Mrowiec}
    \author[T. Rajba]{Teresa Rajba}
    \author[S. Wąsowicz]{Szymon Wąsowicz}
    
    \address{Department of Mathematics, University of Bielsko-Biała, Willowa~2, 43-309 Bielsko-Biała, Poland}

    \email[J. Mrowiec]{jmrowiec@ath.bielsko.pl}
    \email[T. Rajba]{trajba@ath.bielsko.pl}
    \email[S. Wąsowicz]{swasowicz@ath.bielsko.pl}

\begin{abstract}
During the \emph{Conference on Ulam's Type Stability} (Rytro, Po\-land,\linebreak[4] 2014), Ioan Ra\c{s}a recalled his $25$-years-old problem concerning some inequality in\-vol\-ving the Bernstein polynomials. We offer the complete solution (in po\-si\-ti\-ve). As a~tool we use stochastic orderings (which we prove for binomial distributions) as well as so-called concentration inequality. Our methods allow us to pose (and solve) the extended version of the problem in question.
\end{abstract}

\subjclass[2010]{Primary: 26D15; Secondary: 26B25, 60G50} 

\keywords{Bernstein polynomials, stochastic ordering, Ohlin's Lemma, concentration inequality, convex functions, Bernoulli distribution, binomial distribution}
\maketitle
  
\section{Introduction}
The Bernstein fundamental polynomials of degree $n\in\N$ are given by the formula
\[
 b_{n,i}(x)=\binom{n}{i}x^i(1-x)^{n-i}\,,\quad i=0,1,\dots,n\,.
\]
In 2014, during the \emph{Conference on Ulam's Type Stability} held in Rytro (Poland), Ioan Ra\c{s}a recalled his $25$-years-old problem (\cite[Problem~2, p.~164]{CUTS2014}) related to the preservation of convexity by the Bernstein--Schnabl operators.
\begin{pr}
 Prove or disprove that
 \begin{equation}\label{eq:Rasa}
  \sum_{i,j=0}^n\bigl(b_{n,i}(x)b_{n,j}(x)+b_{n,i}(y)b_{n,j}(y)-2b_{n,i}(x)b_{n,j}(y)\bigr)f\left(\frac{i+j}{2n}\right)\ge 0
 \end{equation}
 for each convex function $f\in\C\bigl([0,1]\bigr)$ and for all $x,y\in[0,1]$.
\end{pr}

The aim of this paper is to answer the above-state problem affirmatively (\ie, to prove~\eqref{eq:Rasa}).

Let us invoke some basic notations and results (see \eg~\cite{DenLefSha98}). Let $(\Omega,\mathcal{F},P)$ be a~probability space. As usual, $F_X(x)=P(X<x)$ ($x \in \R$) stands for the probability distribution function of a~random variable~$X:\Omega\to\R$, while $\mu_X$ is the distribution corresponding to~$X$. For real-valued random variables $X,Y$ with finite expectations we say that \emph{$X$ is dominated by~$Y$ in the~stochastic convex ordering sense}, if
\begin{equation}\label{eq:cx}
 \E f(X)\le\E f(Y)
\end{equation}
for all convex functions $f:\R\to\R$ (for which the expectations above exist). In that case we write $X\lcx Y$ or $F_X\lcx F_Y$.

The main idea of our solution is to study the convex stochastic ordering within the class of binomial distributions. To this end we make use of Ohlin's Lemma (\cite[Lemma~2, p. 256]{Ohl69}), which gives a~sufficient condition for two random variables to be in the~stochastic convex ordering relation.
\begin{Ohlin}
 Let $X,Y$ be two random variables and suppose that~$\E X=\E Y$. If the probability distribution functions $F_X,F_Y$ cross exactly once, \ie{},
\begin{equation*}
 F_X(x)\le F_Y(x)\text{ if }x<x_0\quad\text{and}\quad F_X(x)\ge F_Y(x)\text{ if }x>x_0
\end{equation*}
for some $x_0\in\R$, then $X\lcx Y$.
\end{Ohlin}

Originally this lemma was applied to certain insurance problems and it was lesser-known to mathematicians for a long time. It was re-discovered by the se\-cond-named author, who found a~number of applications in the~theory of convex functions (\cf~\cite{Raj14, Raj15}).

\begin{rem}\label{rem:Szostok}
 Szostok noticed in~\cite{Szo15} that if the measures $\mu_X,\mu_Y$ corresponding to $X,Y$, respectively, are concentrated on the interval $[a,b]$, then, in fact, the relation $X\lcx Y$ holds if and only if the inequality~\eqref{eq:cx} is satisfied for all continuous convex functions $f:[a,b]\to\R$.
\end{rem}

Recall that $X\sim B(p)$ means that the random variable $X$ has the Bernoulli distribution with the parameter~$p\in(0,1)$. If $X$ has the binomial distribution with the parameters $n\in\N$ and $p\in(0,1)$ (which we denote by $X\sim B(n,p)$ for short), then, of course,
\begin{equation}\label{eq:Bernnoulli}
 P(X=k)=\binom{n}{k}p^k(1-p)^{n-k}\,,\quad k=0,1,\dots,n\,\quad\text{and}\quad\E X=np.
\end{equation}

Below we recall the binomial convex concentration inequality, which plays an important r\^{o}le in our considerations. It is, in fact, due to Hoeffding~\cite{Hoe63}. Nevertheless, Hoeffding did not state it in the form required for our purposes. The desired form can be found, \emph{e.g.}, in~\cite[Proposition~1, p.~67]{Kle03}.

\begin{thm}\label{th:Hoeffding}
 Let $b_i\sim B(p_i)$ (for $i=1,\dots,n$) be independent random variables. Set $S_n=b_1+\dots+b_n$. Let $\overline{p}=\dfrac{p_1+\dots+p_n}{n}$ and suppose that $S_n^*\sim B(n,\overline{p})$. Then
\begin{equation*}
 \E\Phi(S_n)\le\E\Phi(S_n^*)
\end{equation*}
for any convex function $\Phi:\R\to\R$ (which means that $S_n\lcx S_n^*$).
\end{thm}

A crucial result required to solve Ra\c{s}a's problem reads as follows.
\begin{thm}\label{th:main}
Let $x,y\in(0,1)$ and $n\in\N$. Assume that $X,X_1,X_2, Y,Y_1,Y_2$ are random variables such that $X,X_1,X_2\sim B(n,x)$, $Y,Y_1,Y_2\sim B(n,y)$, $X,Y$ are independent, $X_1,X_2$ are independent and $Y_1,Y_2$ are independent. Then
 \begin{equation}\label{eq:main}
  F_{X+Y}\lcx\frac{1}{2}\bigl(F_{X_1+X_2}+F_{Y_1+Y_2}\bigr)\,.
 \end{equation}
\end{thm}

We postpone the proof to the end of the next section. In Section~\ref{sec2} we present two results on the stochastic convex ordering concerning two binomial distributions. Theorem~\ref{th:Rasa} will be their immediate consequence. Section~\ref{sec3} delivers the solution of the problem of Ra\c{s}a. We note that the inequality~\eqref{eq:main} is no longer valid if we drop the hypothesis that the involved random variables are binomially distributed. In Section~\ref{sec4} we will present the counterexample. In Section~\ref{sec5} we also offer  a~generalization of the Ra\c{s}a problem~\eqref{eq:Rasa} as well as a~generalization of the inequality~\eqref{eq:main} by taking not necessarily two random variables $X,Y$, but the whole family $\bigl\{X_{(1)},\dots,X_{(m)}\bigr\}$ (with $m\ge 2$) of independent random variables.

\section{Stochastic convex ordering --- the case of two binomially distributed random variables}\label{sec2}

This section is devoted to the proof of Theorem~\ref{th:main}. We divide our job into two propositions.

\begin{prop}\label{prop:1}
 Let $x,y\in(0,1)$ and $n\in\N$. Let $X\sim B(n,x)$ and $Y\sim B(n,y)$ be independent random variables. Then
 \[
  X+Y\lcx S_{2n}^*\,,
 \]
 where $S_{2n}^*\sim B\Bigl(2n,\dfrac{x+y}{2}\Bigr).$
\end{prop}
\begin{proof}
 Since $X\sim B(n,x)$ and $Y\sim B(n,y)$ are independent,
 there exist independent random variables $b_1,\dots,b_{2n}$, where $b_i\sim B(p_i)$ for $p_1=\dots=p_n=x$ and $p_{n+1}=\dots=p_{2n}=y$ \st 
 \[
  X=\sum_{i=1}^n b_i\quad\text{and}\quad Y=\sum_{i=n+1}^{2n} b_i\,.
 \]
Set
 \[
  S_{2n}=\sum_{i=1}^{2n}b_i=X+Y\quad\text{and}\quad\overline{p}=\sum_{i=1}^{2n}p_i=\frac{x+y}{2}.
 \]
Now the result follows immediately by Theorem~\ref{th:Hoeffding}.
\end{proof}

\begin{prop}\label{prop:2}
 Let $x,y\in(0,1)$ and $n\in\N$. If $X_1,X_2\sim B(n,x)$ are independent, $Y_1,Y_2\sim B(n,y)$ are independent and $S_{2n}^*\sim B\Bigl(2n,\dfrac{x+y}{2}\Bigr)$, then
 \begin{equation}\label{eq:prop:2}
  F_{S_{2n}^*}\lcx\frac{1}{2}\bigl(F_{X_1+X_2}+F_{Y_1+Y_2}\bigr)\,.
 \end{equation}
\end{prop}
\begin{proof}
 If $x=y$, then $F_{S_{2n}^*}=F_{X_1+X_2}=F_{Y_1+Y_2}$ and~\eqref{eq:prop:2} is trivially satisfied. In the case where $x\ne y$ we assume without loss of generality that $x<y$. Since $X_1,X_2\sim B(n,x)$ are independent as well as $Y_1,Y_2\sim B(n,y)$ are independent, we have $X_1+X_2\sim B(2n,x)$ and $Y_1+Y_2\sim B(2n,y)$. For $k\in\{0,1,\dots,2n\}$ we infer from~\eqref{eq:Bernnoulli} that
 \begin{align*}
  P(X_1+X_2=k)&=\binom{2n}{k}x^k(1-x)^{2n-k}\,,\\[1ex]
  P(Y_1+Y_2=k)&=\binom{2n}{k}y^k(1-y)^{2n-k}\,.
 \end{align*}
 Let us consider the function
 \begin{equation}\label{def:prop2}
  f_{X_1+X_2}(t)=
  \begin{cases}
   \displaystyle\binom{2n}{k}x^k(1-x)^{2n-k}&\text{for }k\le t<k+1,\quad k=0,1,\dots,2n,\\[2ex]
   0&\text{for all other }t.
  \end{cases} 
 \end{equation}
 It is not too difficult to check that
 \begin{equation}\label{eq:prop:2:distr}
  F_{X_1+X_2}(t)=
  \begin{cases}
   0&\text{for }t\le 0\\[1ex]
   \displaystyle\int\limits_0^k f_{X_1+X_2}(u)\dx[u]&\text{for }k-1<t\le k,\quad k=1,2,\dots,2n\,,\\[3ex]
   1&\text{for }t>2n\,.
  \end{cases}
 \end{equation}
 Similarly we define the functions $f_{Y_1+Y_2}$ and $f_{S_{2n}^*}$ by replacing $x$ in the definition of $f_{X_1+X_2}$ by $y,\dfrac{x+y}{2}$, respectively. Of course, a formula analogous to~\eqref{eq:prop:2:distr} holds for the probability distribution functions $F_{Y_1+Y_2}$ and $F_{S_{2n}^*}$.
 
 Now we proceed to the proof of the relation~\eqref{eq:prop:2}. We are going to apply Ohlin's Lemma to the distribution functions $F_{S_{2n}^*}$ and $\dfrac{1}{2}\bigl(F_{X_1+X_2}+F_{Y_1+Y_2}\bigr).$
 Having in mind properties~\eqref{eq:Bernnoulli} we arrive at $\E S_{2n}^*=n(x+y)$ and
 \[
  \dfrac{1}{2}\bigl(\E(X_1+X_2)+\E(Y_1+Y_2)\bigr)=\tfrac{1}{2}(2nx+2ny)=n(x+y)
 \]
 so the distribution functions under consideration admit the same expectations.
 
 The distribution functions $\dfrac{1}{2}\bigl(F_{X_1+X_2}+F_{Y_1+Y_2}\bigr)$ and $F_{S_{2n}^*}$ agree on the interval $(-\infty,0]$ and on the interval $(2n,\infty)$. Then  to verify the second of the hypotheses of Ohlin's Lemma it is enough to prove that there exists $t_0\in(0,2n)$ \st 
 \begin{equation}
  \begin{aligned}\label{ineq:prop2:1}
   \frac{1}{2}\bigl(F_{X_1+X_2}(t)+F_{Y_1+Y_2}(t)\bigr)-F_{S_{2n}^*}(t)&>0 & &\text{for }0<t<t_0\,,\\[1ex]
   \frac{1}{2}\bigl(F_{X_1+X_2}(t)+F_{Y_1+Y_2}(t)\bigr)-F_{S_{2n}^*}(t)&<0 & &\text{for }t_0<t<2n\,.
  \end{aligned}  
 \end{equation}

Since all of the probability distribution functions $F_{X_1+X_2}$, $F_{Y_1+Y_2}$ and $F_{S_{2n}^*}$ are discontinuous at the points $k\in\{0,1,\dots,2n\}$ and constant in between, condition~\eqref{ineq:prop2:1} is equivalent to
 \begin{equation}
  \begin{aligned}\label{ineq:prop2:2}
   \frac{1}{2}\bigl(F_{X_1+X_2}(k)+F_{Y_1+Y_2}(k)\bigr)-F_{S_{2n}^*}(k)&>0 & &\text{for }0<k<t_0\,,\\[1ex]
   \frac{1}{2}\bigl(F_{X_1+X_2}(k)+F_{Y_1+Y_2}(k)\bigr)-F_{S_{2n}^*}(k)&<0 & &\text{for }t_0<k<2n\,.
  \end{aligned}
 \end{equation}
 
Bearing in mind formula~\eqref{eq:prop:2:distr} and the analogous formulae for $F_{Y_1+Y_2}$ and $F_{S_{2n}^*}$ we conclude that the condition~\eqref{ineq:prop2:2} is satisfied if and only if there exist numbers $0<t_1<t_2<2n$ \st
 \begin{equation*}
  \begin{aligned}
   \frac{1}{2}\bigl(f_{X_1+X_2}(k)+f_{Y_1+Y_2}(k)\bigr)-f_{S_{2n}^*}(k)&>0 & &\text{for }0\le k<t_1\,,\\[1ex]
   \frac{1}{2}\bigl(f_{X_1+X_2}(k)+f_{Y_1+Y_2}(k)\bigr)-f_{S_{2n}^*}(k)&<0 & &\text{for }t_1<k<t_2\,,\\[1ex]
   \frac{1}{2}\bigl(f_{X_1+X_2}(k)+f_{Y_1+Y_2}(k)\bigr)-f_{S_{2n}^*}(k)&>0 & &\text{for }t_2<k\le 2n\,.
  \end{aligned}
 \end{equation*}
 
 By~\eqref{def:prop2} and its counterparts for $Y_1+Y_2$ and $S_{2n}^*$, if $k=0,1,\dots,2n$, then
 \begin{multline}\label{113a}
   \frac{1}{2}\bigl(f_{X_1+X_2}(k)+f_{Y_1+Y_2}(k)\bigr)-f_{S_{2n}^*}(k)\\[1ex]
   =\frac{1}{2}\left[\binom{2n}{k}x^k(1-x)^{2n-k}+\binom{2n}{k}y^k(1-y)^{2n-k}\right]\\[1ex]-\binom{2n}{k}\Bigl(\frac{x+y}{2}\Bigr)^k\Bigl(1-\frac{x+y}{2}\Bigr)^{2n-k}=\binom{2n}{k}\psi_k(x,y)\,,
 \end{multline}
 where
 \[
  \psi_k(x,y)=\frac{1}{2}\left[x^k(1-x)^{2n-k}+y^k(1-y)^{2n-k}\right]-\Bigl(\frac{x+y}{2}\Bigr)^k\Bigl(1-\frac{x+y}{2}\Bigr)^{2n-k}\,.
 \]

Consider the case where $k=0$. By strict convexity of the function $u\mapsto(1-u)^{2n}$ on $(0,1)$, we have
\[
  \psi_0(u,v)=\frac{1}{2}\left[(1-u)^{2n}+(1-v)^{2n}\right]-\Bigl(1-\frac{u+v}{2}\Bigr)^{2n}>0
 \]
 for all $u,v\in(0,1)$ with $u\ne v$; in particular 
 \[
  \psi_0(x,y)>0.
 \]
Similarly, for $k=2n$, by the strict convexity of $u\mapsto u^{2n}$ on $(0,1)$, we get
 \[
  \psi_{2n}(u,v)=\frac{1}{2}\left[u^{2n}+v^{2n}\right]-\Bigl(\frac{u+v}{2}\Bigr)^{2n}>0
 \]
for all $u,v\in(0,1)$ with $u\ne v$; in particular
 \[
  \psi_{2n}(x,y)>0.
 \]
Consequently, by \eqref{113a}, 
 \begin{equation*}
  \frac{1}{2}\bigl(f_{X_1+X_2}(k)+f_{Y_1+Y_2}(k)\bigr)-f_{S_{2n}^*}(k)>0
 \end{equation*}
 for $k=0$ and $k=2n$.
 
 Moreover, we \emph{claim} that there exists $k_0\in\{1,2, \dots,2n-1\}$ \st
 \begin{equation}\label{eq:115}
  \frac{1}{2}\bigl(f_{X_1+X_2}(k_0)+f_{Y_1+Y_2}(k_0)\bigr)-f_{S_{2n}^*}(k_0)<0\,.
 \end{equation}
 Assume not. Then
 \[
  \frac{1}{2}\bigl(f_{X_1+X_2}(k)+f_{Y_1+Y_2}(k)\bigr)-f_{S_{2n}^*}(k)>0
 \]
 for all $k\in\{0,1,\dots,2n\}$. Adding these inequalities side by side we arrive at
 \begin{equation*}
  \frac{1}{2}\sum_{k=0}^{2n}f_{X_1+X_2}(k)+\frac{1}{2}\sum_{k=0}^{2n}f_{Y_1+Y_2}(k)>\sum_{k=0}^{2n}f_{S_{2n}^*}(k)\,.
 \end{equation*}
 But using~\eqref{def:prop2} (together with its counterparts for $Y_1+Y_2$ and $S_{2n}^*$) and Newton's Binomial Theorem we get a~contradiction because all the sums above are equal to~$1$. This proves~\eqref{eq:115}.
 
Set $z=\dfrac{x+y}{2}$. By~\eqref{113a}, if $k\in\{0,1,\dots,2n\}$, then
 \begin{equation*}
  \frac{1}{2}\bigl(f_{X_1+X_2}(k)+f_{Y_1+Y_2}(k)\bigr)-f_{S_{2n}^*}(k)=\binom{2n}{k}z^k(1-z)^{2n-k}R_{x,y}(k)\,,
 \end{equation*}
 where
 \[
  R_{x,y}(t)=\frac{1}{2}\Bigl(\frac{1-x}{1-z}\Bigr)^{2n}\left(\frac{\dfrac{x}{1-x}}{\dfrac{z}{1-z}}\right)^t+\frac{1}{2}\Bigl(\frac{1-y}{1-z}\Bigr)^{2n}\left(\frac{\dfrac{y}{1-y}}{\dfrac{z}{1-z}}\right)^t-1\,,\quad t\in[0,2n]\,.
 \]
As $x<y$, we have
 \begin{equation}\label{118}
  \frac{1}{2}\bigl(f_{X_1+X_2}(k)+f_{Y_1+Y_2}(k)\bigr)-f_{S_{2n}^*}(k)>0\iff R_{x,y}(k)>0
 \end{equation}
 for all $k=0,1,\dots,2n$.
\par\noindent
Obviously
 \begin{align}
  R_{x,y}(0)&>0\label{120_1},\\[1ex]
  R_{x,y}(2n)&>0\label{120_2},\\[1ex]
  R_{x,y}(k_0)&<0\quad\text{for some }0<k_0<2n\,.\label{120_3}
 \end{align}
 Since $R_{x,y}^{\prime\prime}(t)>0$ for all $t\in(0,2n)$ (an easy computation), $R_{x,y}$ is a~continuous and convex function on $[0,2n]$. Taking into account \eqref{120_1}, \eqref{120_2} and \eqref{120_3} we conclude that there exist numbers $0<t_1<t_2<2n$ \st 
 \begin{align*}
  R_{x,y}(t)&>0\quad\text{for }0\le t<t_1\,,\\[1ex]
  R_{x,y}(t)&<0\quad\text{for }t_1<t<t_2\,,\\[1ex]
  R_{x,y}(t)&>0\quad\text{for }t_2<t\le 2n\,.
 \end{align*}
 Consequently, by~\eqref{118}, for $k\in\{0, 1, \dots, 2n\}$ we have
 \begin{align*}
  \frac{1}{2}\bigl(f_{X_1+X_2}(k)+f_{Y_1+Y_2}(k)\bigr)-f_{S_{2n}^*}(k)&>0\quad\text{for }0\le k<t_1\,,\\[1ex]
  \frac{1}{2}\bigl(f_{X_1+X_2}(k)+f_{Y_1+Y_2}(k)\bigr)-f_{S_{2n}^*}(k)&<0\quad\text{for }t_1<k<t_2\,,\\[1ex]
  \frac{1}{2}\bigl(f_{X_1+X_2}(k)+f_{Y_1+Y_2}(k)\bigr)-f_{S_{2n}^*}(k)&>0\quad\text{for }t_2<k\le 2n\,.
 \end{align*}
 This implies that the conditions~\eqref{ineq:prop2:2} are satisfied for some $0<t_0<2n$ and the second hypothesis of the Ohlin's Lemma has been verified. Hence~\eqref{eq:prop:2} is satisfied, which completes the proof.\end{proof}

Observe now that Theorem~\ref{th:main} follows immediately from Propositions~\ref{prop:1}~and~\ref{prop:2}.

\section{The problem of Ioan Ra\c{s}a}\label{sec3}

Following Billingsley~\cite{Bil71}, we recall the definition of weak convergence of probability measures. Let~$S$ be a~complete and separable metric space with its Borel $\sigma$-algebra~$\Sigma$. We say that a~sequence $(\mu_m)$ of probability measures on $(S,\Sigma)$ converges weakly to the probability measure $\mu$ (which is denoted by $\mu_m\implies\mu$), if
\begin{equation*}
 \lim_{m\to\infty}\int\limits_S h\dx[\mu_m]=\int\limits_S h\dx[\mu]\quad\text{for all }h\in\C_b(S)\,,
\end{equation*}
where $\C_b(S)$ is the space of all continuous and bounded functions $h:S\to\R$.

For $S=\R$ with the usual topology, if $F_m,F$ are the probability distribution functions of the measures $\mu_m,\mu$, respectively, then $\mu_m\implies\mu$ if and only if $\lim\limits_{m\to\infty}F_m(x)=F(x)$ for all points $x\in\R$ at which~$F$ is continuous.

If $X_m,X:\Omega\to\R$ are random variables ($m\in\N$), then the sequence $(X_m)$ is said to converge weakly to $X$ (write $X_m\implies X$), if the sequence of distributions $(\mu_{X_m})$ converges weakly to~the distribution $\mu_X$ in the above sense.

\begin{rem}\label{Rem31}
 If $\mu_{X_m}$ ($m\in\N$) and $\mu_X$ are concentrated on some compact interval $[a,b]\subset\R$, then $X_m\implies X$ if and only if
 \begin{equation*}
  \lim_{m\to\infty}\int\limits_a^b h\dx[\mu_{X_m}]=\int\limits_a^b h\dx[\mu_X]
 \end{equation*}
 for all continuous functions $h:[a,b]\to\R$.
\end{rem}

A~technical remark will be also needed.

\begin{rem}\label{rem:dzielenie}
 Let $\xi,\eta,\zeta$ be random variables and $a>0$. It is easy to show that
 \begin{align*}
  \xi\lcx\eta&\iff\dfrac{\xi}{a}\lcx\dfrac{\eta}{a}\\
  \intertext{as well as}
  F_{\xi}\lcx\dfrac{1}{2}\bigl(F_{\eta}+F_{\zeta}\bigr)&\iff F_{\frac{\xi}{a}}\lcx\dfrac{1}{2}\bigl(F_{\frac{\eta}{a}}+F_{\frac{\zeta}{a}}\bigr)\,.
 \end{align*}
 Indeed, it is enough to observe that $\E f\left(\frac{\xi}{a}\right)=\E f_{\frac{1}{a}}(\xi)$, where $f_{\frac{1}{a}}(x)=f\left(\frac{x}{a}\right)$ and so on. Of course, $f:\R\to\R$ is convex if and only if $f_{\frac{1}{a}}:\R\to\R$ is convex as well.
\end{rem}

We are now in a~position to achieve the main goal of our paper, which is a~solution of the aforementioned problem of Ioan Ra\c{s}a.
\begin{thm}\label{th:Rasa}
 If $n\in\N$ and
 \[
 b_{n,i}(x)=\binom{n}{i}x^i(1-x)^{n-i}\,,\quad i=0,1,\dots,n\,,
\]
then
 \begin{equation}\label{33}
  \sum_{i,j=0}^n\bigl(b_{n,i}(x)b_{n,j}(x)+b_{n,i}(y)b_{n,j}(y)-2b_{n,i}(x)b_{n,j}(y)\bigr)f\left(\frac{i+j}{2n}\right)\ge 0
 \end{equation}
 for each convex function $f\in\C\bigl([0,1]\bigr)$ and for all $x,y\in[0,1]$.
\end{thm}
\begin{proof}
If $x=y$, then~\eqref{33} is trivially fulfilled, so (by symmetry) it is enough to assume that $0\le x<y\le 1$.

Rewrite~\eqref{33} in the form
\begin{multline*}
 \sum_{k=0}^{2n}\,\sum_{i+j=k} b_{n,i}(x)b_{n,j}(y)f\left(\frac{k}{2n}\right)\\
 \le\frac{1}{2}\sum_{k=0}^{2n}\,\sum_{i+j=k}\bigl(b_{n,i}(x)b_{n,j}(x)+b_{n,i}(y)b_{n,j}(y)\bigr)f\left(\frac{k}{2n}\right)\,,
\end{multline*}
which is equivalent to
\begin{equation}\label{35}
 \E f\left(\frac{X+Y}{2n}\right)\le\frac{1}{2}\left[\E f\left(\frac{X_1+X_2}{2n}\right)+\E f\left(\frac{Y_1+Y_2}{2n}\right)\right]\,,
\end{equation}
where $X_1,X_2$ are independent random variables, $Y_1,Y_2$ are independent random variables and $X,Y$ are independent random variables \st\ four cases are possible:

\begin{enumerate}[(a)]
 \item\label{Rasa_proof_a}
  $0<x<y<1$, $X,X_1,X_2\sim B(n,x)$, $Y,Y_1,Y_2\sim B(n,y)$,
 \item\label{Rasa_proof_b}
  $0=x<y<1$, $\mu_X=\mu_{X_1}=\mu_{X_2}=\delta_0$ ($\delta_{x_0}$ denotes, as usual, the probability measure concentrated at $x_0\in\R$), $Y,Y_1,Y_2\sim B(n,y)$,
 \item\label{Rasa_proof_c}
  $0<x<y=1$, $X,X_1,X_2\sim B(n,x)$, $\mu_Y=\mu_{Y_1}=\mu_{Y_2}=\delta_n$, $\mu_{Y_1+Y_2}=\delta_{2n}$, 
 \item\label{Rasa_proof_d}
  $x=0,y=1$, $\mu_X=\mu_{X_1}=\mu_{X_2}=\delta_0$, $\mu_Y=\mu_{Y_1}=\mu_{Y_2}=\delta_n$.
\end{enumerate}
\par\bigskip
Suppose that~\itemref{Rasa_proof_a} holds. Although we derive
\[
 F_{X+Y}\lcx\frac{1}{2}\bigl(F_{X_1+X_2}+F_{Y_1+Y_2}\bigr)\,,
\]
from Theorem~\ref{th:main}, Remark~\ref{rem:dzielenie} yields
\[
 F_{\frac{X+Y}{2n}}\lcx\frac{1}{2}\left[F_{\frac{X_1+X_2}{2n}}+F_{\frac{Y_1+Y_2}{2n}}\right],
\]
which means that~\eqref{35} holds for all convex functions $f:\R\to\R$, so, by~Remark~\ref{rem:Szostok}, also for all convex functions $f\in\C\bigl([0,1]\bigr)$.
\par\bigskip
Consider now the case~\itemref{Rasa_proof_b}. Let $(x_m)$ be a~sequence of real numbers~\st{} $0<x_m<1$ and $x_m\rightarrow 0$. Let $(\xim)$ be a~sequence of random variables \st{} $\xim\sim B(n,x_m)$ and $\xim,Y$ are independent ($m\in\N$). Let $\xim[(m),1],\,\xim[(m),2]\sim B(n,x_m)$ be independent random variables ($m\in\N$). We shall check that $\mu_{\xim}\implies\delta_0=\mu_X$. Indeed, if $u<0$ then $F_{\xim}(u)=F_X(u)=0$. If $0<u\le 1$ then
\[
 F_{\xim}(u)=P(\xim<u)=P(\xim=0)=(1-x_m)^n\xrightarrow[m\to\infty]{}1.
\]
If $u>1$ then $F_{\xim}(1)\le F_{\xim}(u)\le 1$. Since $F_{\xim}(1)\xrightarrow[m\to\infty]{}1$, then $F_{\xim}(u)\xrightarrow[m\to\infty]{}1$. Hence $\lim\limits_{m\to\infty} F_{\xim}(u)=F_X(u)$ for all $u\ne 0$. Because $F_X$ is continuous at any $u\ne 0$, we get $\mu_{\xim}\implies\delta_0=\mu_X$ (see the introductory note at the beginning of this section). Consequently,
$\mu_{\xim+Y}\implies\mu_{X+Y}$ and $\mu_{\xim[(m),1]+\xim[(m),2]}\implies\mu_{X_1+X_2}$, which implies that
\begin{equation}\label{36}
 \mu_{\frac{\xim+Y}{2n}}\implies\mu_{\frac{X+Y}{2n}}\quad\text{and}\quad\mu_{\frac{\xim[(m),1]+\xim[(m),2]}{2n}}\implies\mu_{\frac{X_1+X_2}{2n}}\,.
\end{equation}
Taking into account $\xim\sim B(n,x_m)$ and $Y\sim B(n,y)$, by the case~\itemref{Rasa_proof_a} we arrive at
\begin{equation}\label{37}
 \E f\left(\frac{\xim+Y}{2n}\right)\le\frac{1}{2}\left[\E f\left(\frac{\xim[(m),1]+\xim[(m),2]}{2n}\right)+\E f\left(\frac{Y_1+Y_2}{2n}\right)\right]
\end{equation}
for all convex functions $f\in\C\bigl([0,1]\bigr)$. Of course, any random variable involved in~\eqref{37} is concentrated on $[0,1]$, so by~\eqref{36} together with Remark~\ref{Rem31} we infer that
\begin{equation}\label{38}
 \begin{aligned}
 \lim\limits_{m\to\infty}\E f\left(\frac{\xim+Y}{2n}\right)&=\E f\left(\frac{X+Y}{2n}\right)\,,\\[1ex]
 \lim\limits_{m\to\infty}\E f\left(\frac{\xim[(m),1]+\xim[(m),2]}{2n}\right)&=\E f\left(\frac{X_1+X_2}{2n}\right)
 \end{aligned}
\end{equation}
for all continuous functions $f:[0,1]\to\R$. The inequality~\eqref{35} follows now by \eqref{37} and~\eqref{38}.

\par\bigskip
In the case~\itemref{Rasa_proof_c} the proof is analogous. Let $(y_m)$ be a~sequence of real numbers~\st~$0<y_m<1$ and $y_m\rightarrow 1$. Let $(\gm)$ be a~sequence of random variables \st~$\gm\sim B(n,y_m)$ and $\gm,X$ are independent. We claim that $\gm\implies Y$. Observe that
\[
 F_{\gm}(n)=P(\gm<n)=1-P(\gm\ge n)=1-P(\gm=n)=1-y_m^n\xrightarrow[m\to\infty]{}0.
\]
For $u<n$, by $0\le F_{\gm}(u)\le F_{\gm}(n)\xrightarrow[m\to\infty]{}0$ we get $F_{\gm}(u)\xrightarrow[m\to\infty]{}0$. If $u>n$ then $F_{\gm}(u)=F_Y(u)=1$. Hence 
$\lim\limits_{m\to\infty} F_{\gm(u)}=F_Y(u)$ for all $u\ne n$. Because $F_Y$ is continuous at any $u\ne n$, we get $\mu_{\gm}\implies\mu_Y$. To prove~\eqref{35} we proceed now similarly as in the case~\itemref{Rasa_proof_b}.

\par\bigskip
Finally, we take into account the case~\itemref{Rasa_proof_d}. The inequality~\eqref{35} could be proved by combining the cases~\itemref{Rasa_proof_b} and~\itemref{Rasa_proof_c}, \ie,~by considering the sequences $(x_m),\,(y_m)$ \st~$x_m,y_m\in(0,1)$, $x_m\rightarrow 0$, $y_m\rightarrow 1$ together with the random variables $\xim\sim B(n,x_m)$ and $\gm\sim B(n,y_m)$. We also notice that
\[
\mu_{\frac{X+Y}{2n}}=\delta_{\frac{1}{2}}\,,\quad
\mu_{\frac{X_1+X_2}{2n}}=\delta_0\quad\text{and}\quad
\mu_{\frac{Y_1+Y_2}{2n}}=\delta_1\,.
\]
Next we could apply Ohlin's Lemma to give an alternative proof, which is considerably easier from the previous one. We omit the details.

\par\bigskip
Thus Theorem~\ref{th:Rasa} is proved and the problem of Ra\c{s}a is completely solved.
\end{proof}

\section{Stochastic convex ordering --- two random variables in a~general case}\label{sec4} 

In this section we show that in the case of any random variables $X,Y$ (not necessarily binomially distributed) the inequality \eqref{eq:main} need not be satisfied. As we can see, Ohlin's Lemma is a strong tool, however, it is worthwhile to notice that in the case of certain inequalities, the corresponding probability distribution functions cross may more than once. Therefore a simple application of Ohlin's Lemma is impossible and an extra idea is needed. To handle such situations, in the papers \cite{OlbSzo15, Szo15}, the authors employed the Levin--Ste\v{c}kin theorem  \cite{LevSte60} (see also \cite{NicPer06}, Theorem 4.2.7).
\begin{LevSte}
Let $a,b\in\R$, $a<b$ and let $F_1,F_2:[a,b]\to\R$ be functions with bounded variation such that $F_1(a)=F_2(a)$. Then, in order that
\begin{equation*}
 \int\limits_a^b f(x)\dx[F_1(x)]\le\int\limits_a^b f(x)\dx[F_2(x)],
\end{equation*}
for all continuous convex functions $f:[a,b]\to\R$, it is necessary and sufficient that $F_1$ and $F_2$ satisfy the following three conditions:
\begin{align*}
 F_1(b)&= F_2(b)\,,\\[1ex]
 \int\limits_a^b F_1(x)\dx&=\int\limits_a^b F_2(x)\dx\,,\\[1ex]
 \int\limits_a^x F_1(t)\dx[t]&\le\int\limits_a^x F_2(t)\dx[t]\quad\text{for all }x\in(a,b)\,.
\end{align*}
\end{LevSte}

To start our considerations, we define the number of sign changes of a function $\varphi:\R\to\R$ by
\[
 S^-(\varphi)=\sup\Bigl\{S^-\bigl[\varphi(x_1),\varphi(x_2),\ldots,\varphi(x_k)\bigr]\,:\,x_1<x_2<\dots< x_k\in\R,\;k\in\N\Bigr\},
\]
where $S^-[y_1,y_2,\ldots,y_k]$ denotes the number of sign changes in the sequence $(y_1,y_2,\ldots,y_k)$ (zero terms are being discarded). Next we say that two real functions $\varphi_1,\varphi_2$ have $n$~crossing points (or cross each other $n$-times) if $S^-(\varphi_1-\varphi_2)=n$. Let $a=x_0<x_1< \dots<x_n<x_{n+1}=b$. The functions  $\varphi_1,\varphi_2$ are said to cross $n$-times at the points $x_1,x_2,\dots,x_n$ (or that $x_1,x_2,\dots,x_n$ are the points of sign changes of $\varphi_1-\varphi_2$) if $S^-(\varphi_1-\varphi_2)=n$ and there exist $a<\xi_1<x_1< \ldots<\xi_n<x_n<\xi_{n+1}<b$ \st~$S^-[\xi_1,\xi_2,\ldots,\xi_{n+1}]=n$.
\par\medskip
The lemma below is due to Szostok (\cf~\cite[Lemma 2]{Szo15}). We quote it in a~sli\-ghtly rewritten form.
\begin{lem}\label{lemma:1.5}
Let $a,b\in\R$, $a<b$ and let $F_1,F_2 \colon [a,b]\to\R$ be  functions with bounded variation such that $F_1(a)=F_2(a)$, $F_1(b)=F_2(b)$, $F=F_2-F_1$, $\int\limits_a^b F(x)\dx=0$. Let $a<x_1<\dots<x_m<b$ be the points of sign changes of $F$. Suppose also that $F(t)\ge 0$ for $t\in(a,x_1)$.
\begin{enumerate}[\upshape(i)]
\item If $m$ is even then the inequality 
 \begin{equation}\label{eq:ord1}
  \int\limits_a^bf(x)dF_1(x)\le\int\limits_a^bf(x)dF_2(x)
\end{equation}
is not satisfied by all continuous convex functions $f:[a,b]\to\R$. 
\item If $m$ is odd, define $A_i$ ($i=0,1,\dots,m$, $x_0=a$, $x_{m+1}=b$) by
\[
 A_i = \int\limits_{x_i}^{x_{i+1}}\left\lvert F(x)\right\rvert\dx.
\]
Then the inequality \eqref{eq:ord1} is satisfied for all continuous convex functions $f:[a,b]\to\R,$  if and only if the following inequalities hold true:
\begin{equation*}
 \begin{split}
  A_0 &\ge A_1\,,\\
  A_0+A_2 & \ge A_1+A_3\,,\\
  &\vdots \\
  A_0+A_2+\dots+A_{m-3} & \ge A_1+A_3+\dots+A_{m-2}\,.
 \end{split}
\end{equation*}
\end{enumerate}
\end{lem}

In a~comment after the statement of~Theorem~\ref{th:main} we indicated that the hypothesis that the random variables involved in the relation~\eqref{eq:main} are binomially distributed is essential. Now we are going to present a counterexample.

\begin{exmp}
 Consider three couples of independent random variables:
 \begin{itemize}
  \item
   $X,Y$ with $\mu_X=\frac{1}{2}(\delta_1+\delta_3)$, $\mu_Y=\frac{1}{2}(\delta_0+\delta_4)$, respectively;
  \item
   $X_1,X_2$ \st~$\mu_{X_1}=\mu_{X_2}=\mu_X$;
  \item
   $Y_1,Y_2$ \st~$\mu_{Y_1}=\mu_{Y_2}=\mu_Y$.
 \end{itemize}
 It is easy to check that
 \begin{align*}
  \mu_{X+Y}&=\tfrac{1}{4}(\delta_1+\delta_3+\delta_5+\delta_7)\,,\\
  \mu_{X_1+X_2}&=\tfrac{1}{4}(\delta_2+\delta_6)+\tfrac{1}{2}\delta_4\,,\\
  \mu_{Y_1+Y_2}&=\tfrac{1}{4}(\delta_0+\delta_8)+\tfrac{1}{2}\delta_4
 \end{align*}
 as well as $\tfrac{1}{2}\bigl(F_{X_1+X_2}+F_{Y_1+Y_2}\bigr)=F_Z$, where $\mu_Z=\tfrac{1}{8}(\delta_0+\delta_2+\delta_6+\delta_8)+\tfrac{1}{2}\delta_4$.

 Put $F=\frac{1}{2}\bigl(F_{X_1+X_2}+F_{Y_1+Y_2}\bigr)-F_{X+Y}$, $a=0$ and $b=8$. Obviously $\int\limits_a^b F(u)\dx[u]=0$. Then $x_1=1$, $x_2=4$ and $x_3=7$ are the points of sign changes of~$F$ and $F(t)\ge 0$ for $t\in(a,x_1)$. Moreover,
 \[
  A_0=\frac{1}{8},\quad A_1=\frac{3}{8},\quad A_2=\frac{3}{8},\quad A_3=\frac{1}{8}\qquad (m=3\text{ is odd}).
 \]
Since $A_0<A_1$, it follows from Lemma~\ref{lemma:1.5} that the relation~\eqref{eq:main}, \ie
 \[
  F_{X+Y}\lcx\frac{1}{2}\bigl(F_{X_1+X_2}+F_{Y_1+Y_2}\bigr)\,,
 \]
 does not hold.
\end{exmp}

\section{An extension of the problem of Ra\c{s}a}\label{sec5}

Let us start with the extension of the results of Propositions~\ref{prop:1}~and~\ref{prop:2} to the case of any~finite number of independent random variables.

\begin{prop}\label{prop:3}
 Let $m,n\in\N$, $m\ge 2$, $x_1,\dots,x_m\in(0,1)$. Suppose that
 \begin{enumerate}[\upshape(i)]
  \item
   $X_{(1)},\dots,X_{(m)}$ are independent random variables \st~$X_{(i)}\sim B(n,x_i)$, $i=1,\dots,n$;
  \item
   $S_{mn}^*\sim B(mn,\bar{x})$, where $\bar{x}=\frac{1}{m}\displaystyle\sum_{i=1}^m x_i$;
   \item
    $X_{(i),1},\dots,X_{(i),m}$ are independent random variables \st~$X_{(i),j}\sim B(n,x_i)$, $j=1,\dots,m$, $i=1,\dots,n$.
 \end{enumerate}
 Then
 \begin{align}
  X_{(1)}+\dots+X_{(m)}&\lcx S_{mn}^*,\label{prop:3_a}\\[1ex]
  F_{S_{mn}^*}&\lcx\frac{1}{m}\Bigl[F_{X_{(1),1}+\dots+ X_{(1),m}} +\dots+F_{X_{(m),1}+\dots+X_{(m),m}}\Bigr],\label{prop:3_b}\\[1ex]
 F_{ X_{(1)}+\dots+X_{(m)}}&\lcx\frac{1}{m}\Bigl[F_{X_{(1),1}+\dots+ X_{(1),m}} +\dots+F_{X_{(m),1}+\dots+X_{(m),m}}\Bigr].\label{prop:3_c}
 \end{align}
\end{prop}

\begin{proof}
 It is enough to prove the relations \eqref{prop:3_a}, \eqref{prop:3_b}. Relation~\eqref{prop:3_c} is their immediate consequence.
 \par\bigskip
 The proof of relation~\eqref{prop:3_a} is short. Assuming that all the hypotheses are satisfied, there exist independent random variables $b_1,b_2,\dots,b_{mn}$ \st
 \[
  b_j\sim B(x_i)\,,\qquad i=1,\dots,m,\quad j=(i-1)n+1,\dots,in
 \]
 and
 \[
  X_{(i)}=\sum_{j=(i-1)n+1}^{in} b_j\,,\qquad i=1,\dots,m\,.
 \] 
  Then $X_{(1)}+\dots+X_{(m)}=\displaystyle\sum_{j=1}^{mn}b_j$ and~\eqref{prop:3_a} is an immediate consequence of Theorem~\ref{th:Hoeffding}.
  \par\medskip
  The rest of the proof is devoted to the relation~\eqref{prop:3_b}.
  \par\medskip
  Because if $x_1=\dots=x_m$ then~\eqref{prop:3_b} is trivially satisfied, assume that this condition does not hold. Without loss of generality assume moreover that

 \begin{equation*}
  x_1\le x_2\dots\le x_m\quad\text{and }x_1<x_m\,.
 \end{equation*}
 Let $i=1,2,\dots,m$. Since $X_{(i)1},X_{(i)2},\dots,X_{(i)m}\sim B(n,x_i)$ are independent, we have $X_{(i)1}+X_{(i)2}+\dots+X_{(i)m}\sim B(mn,x_i)$. Hence
 \[
  P\bigl(X_{(i)1}+X_{(i)2}+\dots+X_{(i)m}=k\bigr)=\binom{mn}{k}x_i^k(1-x_i)^{mn-k}\,,\qquad k=0,1,\dots,mn.
 \]
 For the function
 \begin{multline}\label{53}
  f_{X_{(i)1}+X_{(i)2}+\dots+X_{(i)m}}(t)\\
  =
  \begin{cases}
   \displaystyle\binom{mn}{k}x_i^k(1-x_i)^{mn-k}&\text{for }k\le t<k+1,\;k=0,1,\dots,mn\,,\\[2ex]
   0&\text{for all other }t
  \end{cases}
 \end{multline}
 we easily check that
 \begin{multline}\label{54}
  F_{X_{(i)1}+X_{(i)2}+\dots+X_{(i)m}}(t)\\
  =
  \begin{cases}
   0&\text{for }t\le 0\,,\\[1ex]
   \displaystyle\int\limits_0^k f_{X_{(i)1}+X_{(i)2}+\dots+X_{(i)m}}(u)\dx[u]&\text{for }k-1< t\le k,\;k=1,\dots,mn\,,\\[2ex]
   1&\text{for }t>mn\,.
  \end{cases}
 \end{multline}
 Similarly we define the function $f_{S_{mn^*}}$ by putting $\bar{x}$ instead of~$x_i$ in the definition~\eqref{53}. Of course, the formula analogous to~\eqref{54} holds for the distribution function~$F_{S_{mn^*}}$.
 \par\medskip
 As in the proof of~Proposition~\ref{prop:2} (\ie{} in the case $m=2$), now we check the hypotheses of Ohlin's Lemma. The first one (concerning the equality of expectations) is easily fulfilled, so we turn our attention to the second one. It is enough to prove that there exists $t_0\in(0,mn)$ \st 
 \begin{equation}
  \begin{aligned}\label{55}
   \frac{1}{m}\Bigl[F_{X_{(1),1}+\dots+ X_{(1),m}}(k) +\dots+F_{X_{(m),1}+\dots+X_{(m),m}}(k)\Bigr]-F_{S_{mn}^*}(k)&>0\\ 
    \text{for }&0<k<t_0\,,\\[1ex]
   \frac{1}{m}\Bigl[F_{X_{(1),1}+\dots+ X_{(1),m}}(k) +\dots+F_{X_{(m),1}+\dots+X_{(m),m}}(k)\Bigr]-F_{S_{mn}^*}(k))&<0\\ 
    \text{for }&t_0<k<mn\,.
  \end{aligned}  
 \end{equation}
 Having in mind~\eqref{54} and the analogous formula for $F_{S_{mn}^*}$ we infer that condition~\eqref{55} is satisfied if and only if there exist $0<t_1<t_2<mn$ \st{}
 \begin{equation}
  \begin{aligned}\label{56}
   \frac{1}{m}\Bigl[f_{X_{(1),1}+\dots+ X_{(1),m}}(k) +\dots+f_{X_{(m),1}+\dots+X_{(m),m}}(k)\Bigr]-f_{S_{mn}^*}(k)&>0\\ 
    \text{for }&0\le k<t_1\,,\\[1ex]
   \frac{1}{m}\Bigl[f_{X_{(1),1}+\dots+ X_{(1),m}}(k) +\dots+f_{X_{(m),1}+\dots+X_{(m),m}}(k)\Bigr]-f_{S_{mn}^*}(k)&<0\\ 
    \text{for }&t_1<k<t_2\,,\\[1ex]
   \frac{1}{m}\Bigl[f_{X_{(1),1}+\dots+ X_{(1),m}}(k) +\dots+f_{X_{(m),1}+\dots+X_{(m),m}}(k)\Bigr]-f_{S_{mn}^*}(k)&>0\\ 
    \text{for }&t_2<k\le mn\,. 
  \end{aligned}  
 \end{equation}
 By~\eqref{53} and its counterpart for $f_{S_{mn}^*}$, if $k=0,1,\dots,mn$ then
 \begin{multline}\label{57}
  \frac{1}{m}\Bigl[f_{X_{(1),1}+\dots+ X_{(1),m}}(k) +\dots+f_{X_{(m),1}+\dots+X_{(m),m}}(k)\Bigr]-f_{S_{mn}^*}(k)\\
  =\frac{1}{m}\left[\binom{mn}{k}x_1^k(1-x_1)^{mn-k}+\dots+\binom{mn}{k}x_m^k(1-x_m)^{mn-k}\right]\\
  -\binom{mn}{k}\bar{x}^k(1-\bar{x})^{mn-k}=\binom{mn}{k}\psi_k(x_1,\dots,x_m)\,,
 \end{multline}
 where
 \begin{multline}\label{eq:psi_k}
  \psi_k(x_1,\dots,x_m)\\=\frac{1}{m}\left[x_1^k(1-x_1)^{mn-k}+\dots+x_m^k(1-x_m)^{mn-k}\right]-\bar{x}^k(1-\bar{x})^{mn-k}.
 \end{multline}
 If $k=0$ or $k=mn$ then $\psi_k(x_1,\dots,x_m)>0$ by the strict convexity (on $(0,1)$) of the functions $u\mapsto (1-u)^{mn}$ and $u\mapsto u^{mn}$, respectively. By~\eqref{57} we get
 \[
  \frac{1}{m}\Bigl[f_{X_{(1),1}+\dots+ X_{(1),m}}(k) +\dots+f_{X_{(m),1}+\dots+X_{(m),m}}(k)\Bigr]-f_{S_{mn}^*}(k)>0\,.
 \]
 Similarly as in the proof of Proposition~\ref{prop:2} we can show that there exists $k_0\in\{1,2,\dots,mn-1\}$ \st{}
 \[
  \frac{1}{m}\Bigl[f_{X_{(1),1}+\dots+ X_{(1),m}}(k_0) +\dots+f_{X_{(m),1}+\dots+X_{(m),m}}(k_0)\Bigr]-f_{S_{mn}^*}(k_0)<0\,.
 \]
Taking into account~\eqref{eq:psi_k}, for $k=0,1,\dots,mn$ we have
\begin{equation}\label{58}
  \psi_k(x_1,\dots,x_m)=\bar{x}^k(1-\bar{x})^{mn-k}R_{x_1,\dots,x_m}(k)\,,
 \end{equation}
 where
 \[
  R_{x_1,\dots,x_m}(t)=\frac{1}{m}\left[\left(\frac{1-x_1}{1-\bar{x}}\right)^{mn}\left(\dfrac{\dfrac{x_1}{1-x_1}}{\dfrac{\bar{x}}{1-\bar{x}}}\right)^t+\dots+\left(\frac{1-x_m}{1-\bar{x}}\right)^{mn}\left(\dfrac{\dfrac{x_m}{1-x_m}}{\dfrac{\bar{x}}{1-\bar{x}}}\right)^t\right]-1
 \]
for $t\in[0,mn]$. 
By computing the second derivative we convince ourselves that this function is strictly convex on $(0,mn)$. Then $R_{x_1,\dots,x_m}(0)>0$, $R_{x_1,\dots,x_m}(mn)>0$ and $R_{x_1,\dots,x_m}(k_0)<0$ for some $0<k_0<mn$. Combining this with continuity of $R_{x_1,\dots,x_m}$ on $[0,mn]$ we conclude that there exist $0<t_1<t_2<mn$ \st
 \begin{align*}
  R_{x_1,\dots,x_m}(t)&>0\quad\text{for }0\le t<t_1\,,\\
  R_{x_1,\dots,x_m}(t)&<0\quad\text{for }t_1<t<t_2\,,\\
  R_{x_1,\dots,x_m}(t)&>0\quad\text{for }t_2<t\le mn\,.
 \end{align*}
Following \eqref{57} and \eqref{58} we see that the relations \eqref{56} hold, so the second hypothesis of Ohlin's Lemma is fulfilled. It is now enough to apply this result to complete the proof of Proposition~\ref{prop:3}.
\end{proof}

Now we present the result which extends Theorem~\ref{th:Rasa}, and, therefore, generalizes the problem of~Ra\c{s}a.

\begin{thm}
 Let $m,n\in\N$, $m\ge 2$, $x_1,\dots,x_m\in(0,1)$. Then
 \begin{multline}\label{514}
  \sum_{i_1,\dots,i_m=0}^n\bigl(b_{n,i_1}(x_1)\cdots b_{n,i_m}(x_1)+\dots+b_{n,i_1}(x_m)\dots b_{n,i_m}(x_m)\\
  -mb_{n,i_1}(x_1)\dots b_{n,i_m}(x_m)\bigr)f\left(\frac{i_1+\dots+i_m}{mn}\right)\ge 0
 \end{multline}
 for each convex function $f\in\C\bigl([0,1]\bigr)$ and for all $x_1,\dots,x_m\in[0,1]$.
\end{thm}

\begin{proof}
 Rewrite~\eqref{514} in the form
 \begin{align*}
  &\sum_{k=0}^{mn}\;\sum_{i_1+\dots+i_m=k}b_{n,i_1}(x_1)\cdots b_{n,i_m}(x_m)f\left(\frac{k}{mn}\right)\\
  \le&\frac{1}{m}\sum_{k=0}^{mn}\;\sum_{i_1+\dots+i_m=k}\bigl(b_{n,i_1}(x_1)\cdots b_{n,i_m}(x_1)+\dots+b_{n,i_1}(x_m)\dots b_{n,i_m}(x_m)\bigr)f\left(\frac{k}{mn}\right)\,,
 \end{align*}
 which is equivalent to
 \begin{multline}\label{516}
  \E f\left(\frac{X_{(1)}+\dots+X_{(m)}}{mn}\right)\\
  \le\frac{1}{m}\left[\E f\left(\frac{X_{(1),1}+\dots+X_{(1),m}}{mn}\right)+\dots+\E f\left(\frac{X_{(m),1}+\dots+X_{(m),m}}{mn}\right)\right]\,,
 \end{multline}
 where $X_{(1)},\dots,X_{(m)}$ are independent random variables and for all $i\in\{1,\dots,m\}$ the random variables $X_{(i),1},\dots,X_{(i),m}$ are independent with
 \begin{align*}
  X_{(i)},X_{(i),1},\dots,X_{(i),m}\sim B(n,x_i)\,,&\quad\text{if }x_i\in(0,1)\,,\\
  \mu_{X_{(i)}}=\mu_{X_{(i),1}}=\dots=\mu_{X_{(i),m}}=\delta_0\,,&\quad\text{if }x_i=0\,,\\
  \mu_{X_{(i)}}=\mu_{X_{(i),1}}=\dots=\mu_{X_{(i),m}}=\delta_n\,,&\quad\text{if }x_i=1\,.
 \end{align*}

The proof of the inequality~\eqref{516} (based on Proposition~\ref{prop:3}) is similar to the proof of Proposition~\ref{prop:2} (\ie{} in the case $m=2$) and we omit it.
\end{proof}

\bibliographystyle{plain}

\begin{thebibliography}{10}

\bibitem{CUTS2014}
Report of meeting, {C}onference on {U}lam’s {T}ype {S}tability {R}ytro,
  {P}oland, {J}une 2-6, 2014.
\newblock {\em Ann. Univ. Paedagog. Crac. Stud. Math.}, 13:139--169, 2014.

\bibitem{Bil71}
Patrick Billingsley.
\newblock {\em Weak convergence of measures: {A}pplications in probability}.
\newblock Society for Industrial and Applied Mathematics, Philadelphia, Pa.,
  1971.
\newblock Conference Board of the Mathematical Sciences Regional Conference
  Series in Applied Mathematics, No. 5.

\bibitem{DenLefSha98}
Michel {Denuit}, Claude {Lefevre}, and Moshe {Shaked}.
\newblock {The $s$-convex orders among real random variables, with
  applications.}
\newblock {\em {Math. Inequal. Appl.}}, 1(4):585--613, 1998.

\bibitem{Hoe63}
W.~{Hoeffding}.
\newblock {Probability inequalities for sums of bounded random variables.}
\newblock {\em {J. Am. Stat. Assoc.}}, 58:13--30, 1963.

\bibitem{Kle03}
Thierry Klein.
\newblock In\'egalit\'es de concentration, martingales et arbres al\'eatoires.
\newblock {\em Th\`ese, Universit\'e de Versailles-Saint-Quentin}, 2003.

\bibitem{LevSte60}
V.~I. Levin and S.~B. Ste{\v{c}}kin.
\newblock Inequalities.
\newblock {\em Amer. Math. Soc. Transl. (2)}, 14:1--29, 1960.

\bibitem{NicPer06}
Constantin~P. Niculescu and Lars-Erik Persson.
\newblock {\em Convex functions and their applications}.
\newblock CMS Books in Mathematics/Ouvrages de Math\'ematiques de la SMC, 23.
  Springer, New York, 2006.
\newblock A contemporary approach.

\bibitem{Ohl69}
Jan Ohlin.
\newblock On a class of measures of dispersion with application to optimal
  reinsurance.
\newblock {\em ASTIN Bulletin}, 5:249--266, 1969.

\bibitem{OlbSzo15}
Andrzej Olbry{\'s} and Tomasz Szostok.
\newblock Inequalities of the {H}ermite-{H}adamard type involving numerical
  differentiation formulas.
\newblock {\em Results Math.}, 67(3-4):403--416, 2015.

\bibitem{Raj14}
Teresa Rajba.
\newblock On the {O}hlin lemma for {H}ermite-{H}adamard-{F}ej\'er type
  inequalities.
\newblock {\em Math. Inequal. Appl.}, 17(2):557--571, 2014.

\bibitem{Raj15}
Teresa Rajba.
\newblock On strong delta-convexity and {H}ermite-{H}adamard type inequalities
  for delta-convex functions of higher order.
\newblock {\em Math. Inequal. Appl.}, 18(1):267--293, 2015.

\bibitem{Szo15}
Tomasz Szostok.
\newblock Ohlin's {L}emma and some inequalities of the {H}ermite-{H}adamard
  type.
\newblock {\em Aequationes Math.}, 89:915--926, 2015.

\end{thebibliography}

\end{document}